\documentclass{amsart}
\usepackage{amsmath, amssymb, verbatim, hyperref}

\def\Ric{\mathop{\rm Ric}}
\def\cRic{\mathop{\rm R\i\makebox[0pt]{\raisebox{5pt}{\tiny$\circ$\;\,}}c}}

\def\dist{\mathop{\rm dist}}

\def\Riem{\mathop{\rm Rm}}

\def\Vol{\mathop{\rm Vol}}

\def\supp{\mathop{\rm supp}}

\def\PPP{\mathop{\mathcal{P}}}
\def\RRR{\mathop{\mathcal{P}}}

\def\be{\begin{eqnarray}}
\def\ee{\end{eqnarray}}
\def\beg{\begin{eqnarray*}}
\def\ees{\end{eqnarray*}}

\def\XXint#1#2#3{{\setbox0=\hbox{$#1{#2#3}{\int}$ }
\vcenter{\hbox{$#2#3$ }}\kern-.5\wd0}}

\newtheorem{theorem}{Theorem}[section]
\newtheorem{proposition}[theorem]{Proposition}
\newtheorem{lemma}[theorem]{Lemma}

\newcounter{myremark}[section]

\date{\small\it August 22, 2017}
\author{Brian Weber}
\title{Energy and Asymptotics of Ricci-Flat 4-Manifolds with a Killing Field}


\begin{document}

\maketitle

\begin{abstract}
	If a complete 4-manifold with $\Ric=0$ has a nowhere zero Killing field, we prove it is flat, generalizing a classic result on compact manifolds.
	If the Killing field has compact zero-locus, we compute the manifold's $L^2$ energy.
\end{abstract}

\section{Introduction}

We prove that a complete Ricci-flat 4-manifold with a nowhere zero Killing field is flat, providing an affirmative answer to a conjecture stated in \cite{And2}.
If the zero-set of the Killing field is compact, we compute the $L^2$ curvature energy of the manifold in terms of its asymptotic volume ratio.

\begin{theorem} \label{ThmI}
Assume $(N^4,\,g)$ is a complete, Ricci-flat 4-manifold with a nowhere-zero Killing field $X$.
Then $N$ is flat.
\end{theorem}
\begin{theorem} \label{ThmIII}
Assume $(N^4,g)$ is complete, Ricci-flat, and has a Killing field $X$.
Then
\begin{eqnarray}
\frac{1}{8\pi^2}\int_{N^4}|\Riem|^2=\epsilon+\sum_{i=1}^m\chi({N}_{X,i}) \,-\,
\lim_{s\rightarrow\infty}\frac{Vol\,B(p,s)}{\frac12\pi^2s^4}. \label{IneqMainRiemEst2}
\end{eqnarray}
where $\{N_{X,i}\}_{i=1}^m$ are the components of the null-set of $X$, and $\epsilon$ is defined as follows: if $N^4$ is not ALE then $\epsilon=0$, and if $N^4$ is ALE, then $\epsilon=-1+Ind_{\infty}X$ where $Ind_{\infty}X$ is the index of $X$ at infinity.
\end{theorem}

These results have intrinsic interest, as they fill in a gap in the understanding of Ricci-flat manifolds.
But we are motivated by the use of Theorem \ref{ThmI} in ruling out certain ``blow-up'' instantons.
If an instanton arises as some cover of a region of collapse with locally bounded curvature, then a nowhere zero Killing fields automatically appears in the limit; if this limit has zero Ricci curvature, it is flat. In this way collapsing with bounded Ricci curvature but unbounded sectional curvature can be ruled out.
See \cite{Web2} for the use of Theorem \ref{ThmI} in establishing a regularity estimate.

Before outlining the proofs, we define the {\it s-local curvature radius}
\begin{eqnarray*}
r_{\mathcal{R}}^s(p)\;=\;\sup\,\{\,0<r<s\;\big|\;|\Riem|<r^{-2}\;{\rm on}\;{B}(p,r)\;\},
\end{eqnarray*}
and $r_{\mathcal{R}}=r_{\mathcal{R}}^\infty$; compare (1.7) and (2.2) of \cite{CT}.
Comparing the expressions for $|\Riem|^2$ and the Chern-Gauss-Bonnet 4-form $\mathcal{P}_\chi$
\begin{eqnarray*}
	\begin{aligned}
		|\Riem|^2&\;=\;\frac16R^2\,+\,2|\cRic|^2\,+\,|W|^2, \quad and \\
		8\pi^2\mathcal{P}_\chi&\;=\;\left(\frac{1}{24}R^2\,-\,\frac12|\cRic|^2\,+\,\frac14|W|^2\right)\,dVol,
	\end{aligned}
\end{eqnarray*}
we arrive at the classic observation that $\cRic=0$ implies $\int|\Riem|^2dVol$ and $\int\mathcal{P}_\chi$ are equivalent norms.
The Killing field $X$ can be used to explicitly transgress the canonical Chern-Gauss-Bonnet 4-form $\mathcal{P}_\chi$, and therefore the curvature 4-form: $|\Riem|^2dVol=d\,\mathcal{TP}_\chi$, where the 3-form $\mathcal{TP}_\chi$ is a combination of $X$, $\nabla{X}$, and curvature terms.
Approximately $|\mathcal{TP}_\chi|=|W|\frac{|\nabla{X}|}{|X|}+\frac{|\nabla{X}|^3}{|X|^3}$; the exact expression is (\ref{EqnEulerTransgression}).
Integration by parts gives $\int\varphi|\Riem|^2dol=\sum\chi(N_{X,i})-\int{d}\varphi\wedge\mathcal{TP}$, provided $\{\varphi=1\}$ contains the zero locus of $X$.

Then we estimate $\int{}d\varphi\wedge\mathcal{TP}_\chi$ back in terms of $\int|\Riem|^2$.
As we show in section \ref{SecTransgressions},
\begin{eqnarray}
	\int{d}\varphi\wedge\mathcal{TP}_\chi \;\le\;
	\int|d\varphi|(r_{\mathcal{R}}{})^{-3}dVol. \label{IneqInitialTransBound}
\end{eqnarray}
This is not straighforward, however, as the value $|\mathcal{TP}_\chi|$ depends on $|\nabla{X}|/|X|$, over which we do not have any pointwise control.
But although a pointwise estimate is impossible, an integral estimate is not; this the subject of \S\ref{SubSecFlatCase}.
To obtain (\ref{IneqInitialTransBound}), we divide the manifold into regions where $|\nabla{X}|/|X|$ is not too large so a pointwise estimate suffices, and regions where $|\nabla{X}|/|X|$ is very large but where (as we show) we have enough structural information on $X$  to estimate the necessary integral.

Having obtained (\ref{IneqInitialTransBound}), we use the Cheeger-Tian method from \cite{CT}.
The function $(r_{\mathcal{R}})^{-1}$ is not much different from a Hardy-Littlewood maximal function for $|\Riem|^2$, and a weak-(1,1) inequality lets us estimate the right side of (\ref{IneqInitialTransBound}) in terms of $\left(\int|\Riem|^2\right)^{3/4}$---note the improved exponent.
An iteration argument similar to that of \cite{CT} then proceeds, and we find that $\int|\Riem|^2$ is uniformly bounded in terms of $\sum_i\chi(N_{X,i})$ and a definite multiple of the manifold's asymptotic volume ratio.
In the case that the manifold has zero asymptotic volume ratio, the proof is then complete.

The remaining case, explored in \S\ref{SubSecALECase}, is that the manifold is ALE.
The orbits of $X$ must be bounded, or else necessarily $\int|W|^2$ is zero or infinite, therefore zero because we already obained a universal estimate of $\int|W|^2$ in terms of the asymptotic volume ratio.
If the orbits are bounded, then the Euler number of its universal cover is zero, which we show, with a simple homological argument, is also impossible.

\stepcounter{myremark} 
{\bf Remark \arabic{section}.\arabic{myremark}}. Theorem \ref{ThmI} was conjectured in \cite{And2} (remark 4.3) and \cite{And4}.
The rationale seemed to be via analogy with the Lorenzian case, where geodesically complete spacetime instantons with a time-like Killing field are flat (theorem 0.1 of \cite{And3}).

\stepcounter{myremark}
{\bf Remark \arabic{section}.\arabic{myremark}}. 
Under the strong additional assumption that the manifold is ALE (asymptotically locally Euclidean) or ALF (asymptotically locally flat), the equivalent of Theorem \ref{ThmIII} has long been known: see, for instance, equation (12) in \cite{GPR} or the discussion of characteristic classes and instantons in \cite{EGH}.
What we have done is remove the strong hypotheses on the asymptotic geometry of the manifold.

\stepcounter{myremark}
{\bf Remark \arabic{section}.\arabic{myremark}}.
As indicated above, the necesssary computations are complicated by the fact that the quantity $|\nabla{X}|/|X|$ has no apriori control.
But in some important cases this quantity {\it can} be controlled, for instance if $X$ represents an F-structure of rank 1 that comes from collapse with locally bounded curvature \cite{CG1} \cite{CT}.
A straighforward criterion for simplifying the computations would be to require $|\nabla{X}|/|X|\le{}Ar_{\mathcal{R}}(p)^{-1}$ where $A$ is some constant; indeed this is met in the case of F-structures arising from collapse with locally bounded curvature \cite{CT}.
If this condition holds, then the inequality (\ref{IneqInitialTransBound}) is immediate, and \S\ref{SubSecFlatCase} can be skipped.

\stepcounter{myremark}
{\bf Remark \arabic{section}.\arabic{myremark}}. 
The methods of this paper suffice to weaken the hypothesis that $\Ric=0$ in Theorem \ref{ThmI} to only requiring that $\Ric\ge0$ and the trace-free Ricci be sufficiently pinched: $|\cRic|^2\le\frac{1}{24}{R}^2$.
In fact even $\Ric\ge0$ be relaxed by requiring that the $L^p$-norm of $\Ric_-$ be controlled for some $p>2$.
For the sake of clarity, we do not pursue this level of generality.

{\bf Acknowledgement}. The author wishes to thank Mike Anderson for bringing the conjecture that became Theorem \ref{ThmI} to his attention.

\section{Computing the Transgression} \label{SecTransgressions}

This section's purpose is twofold: we explicitly compute the transgression $\mathcal{TP}_\chi$ given a Killing field, and we evaluate the integral of this transgression in a definite region where the ratio $|\nabla{X}|/|X|$ is extremely large but not infinite.

\subsection{General construction of transgressions}

Given a vector field, a standard process allows any characteristic 4-form to be transgressed in an essentially canonical manner.
When the vector field is Killing, the resulting expressions are simpler than in the usual textbook cases.
Let $X$ be a Killing field on an $n$-manifold with dual 1-form $X_\flat$.
Set $K\;=\;|X|^{-2}\,X_\flat\otimes\nabla{X}$; because $\left<\nabla{X},\,\cdot\,\right>$ is antisymmetric, $K$ is an $\mathfrak{so}(n)$-valued 1-form, so $\widetilde\nabla=\nabla-K$ is a new connection.
Also $X$ is covariant-constant in $\widetilde\nabla$, as $\widetilde\nabla{X}
=\nabla{X}-|X|^{-2}X_\flat(X)\nabla{X}=0.$

Let $\{e_1=X|X|^{-1},\,e_2,\,\dots,\,e_n\}$ be a local frame with $[X,\,e_i]=0$.
With $\mathcal{L}$ denoting the Lie derivative, the connection 1-forms $\widetilde{A}=A-K$ then obey $\mathcal{L}_X\widetilde{A}=\mathcal{L}_X\widetilde\nabla=0$ and $i_X\widetilde{A}=0$, so also $i_X[\widetilde{A},\widetilde{A}]=0$.
Thus the associated curvature 2-forms have a null vector:
\begin{eqnarray*}
	\begin{aligned}
		i_X\widetilde{F}
		&\;=\;i_Xd\widetilde{A}\,+\,\frac12i_X[\widetilde{A},\widetilde{A}]
		\;=\;\mathcal{L}_X\widetilde{A}\;=\;0.
	\end{aligned}
\end{eqnarray*}
Letting $\PPP$ be any symmetric invariant polynomial in $\mathfrak{so}(n)$ that is homogeneous of order $2$, we know that $\PPP(\underline{F},\underline{F})$ represents a top characteristic class on $N^4$, where $\underline{F}$ is the curvature of any connection.
Because $\widetilde{F}$ has a null vector, $\mathcal{P}(\widetilde{F},\widetilde{F})=0$.

Put $A_t=A-tK$ and $F_t=dA_t+\frac12[A_t,\,A_t]$.
Since $K$ has the form $K=v_\flat\otimes{h}$ for $v_\flat=|X|^{-2}X_\flat$ and some $h\in\mathfrak{so}(n)$, a consequence is that $[K,K]=0$.
Using this, we see $D_tK=DK$ and $F_t=F\,-\,t\,DK$.
The standard computation gives $\frac{d}{dt}\mathcal{P}(F_t,F_t)
=-2\mathcal{P}(D_tK,F_t)
=-2d\mathcal{P}(K,F_t)$.
Therefore
$\PPP(F,\,F)\,-\,\PPP(\widetilde{F},\,\widetilde{F}) \;=\;2d\int_0^1\mathcal{P}(K,\,F_t)\,dt.$
Therefore if we set $\mathcal{TP}=2\int_0^1\mathcal{P}(K,F_t)dt$, we have $\mathcal{P}=d\,\mathcal{TP}$.
Explicitly,
\begin{eqnarray}
\mathcal{TP}\;=\;2\mathcal{P}(K,\,F) \,-\,\mathcal{P}(K,\,DK). \label{EqnTPDim4}
\end{eqnarray}

\subsection{The Chern-Gauss-Bonnet transgression form}

We must explicitly compute $\nabla{K}$ and its antisymmetrization $DK\in\bigwedge^2\otimes\mathfrak{so}(4)$.
Letting $A$, $B$, and $w$ be vector fields, we have
\begin{eqnarray*}
	\begin{aligned}
		&(\nabla{K})(A,B).w \;=\;\nabla_A\left(K(B).w\right)\,-\,K(\nabla_AB).w\,-\,K(B).\nabla_Aw\\
		&\quad=2|X|^{-4}\left<\nabla_XX,A\right>\left<X,B\right>\nabla_wX
		+|X|^{-2}\left<\nabla_AX,B\right>\nabla_wX +|X|^{-2}\left<X,B\right>\nabla^2_{A,w}X.
	\end{aligned}
\end{eqnarray*}
Since $dX_\flat(A,B)=2\left<\nabla_AX,\,B\right>$, antisymmetrizing in the $A$, $B$ positions gives
\begin{eqnarray}
	DK\;=\;\left(|X|^{-4}\left(i_XdX_\flat\right)\wedge{X}_\flat
	\,+\,|X|^{-2}dX_\flat\right)\otimes\nabla{X} \,-\,|X|^{-2}X_\flat\wedge\nabla^2X
\end{eqnarray}
where $X_\flat\wedge{}\nabla^2X$ means $X^sg_{si}X^l{}_{,kj}dx^i\wedge{}dx^j\otimes{}dx^k\otimes{}\frac{\partial}{\partial{x}^l}\in\bigwedge{}^2\otimes\mathfrak{so}(4)$.
Using the Hodge star, $*:\mathfrak{so}(4)\rightarrow\mathfrak{so}(4)$, and the fact that $\Riem=W$ (since $\Ric=0$), we have
\be
\begin{aligned}
&\quad{}K\wedge*{D}K
\;=\;|X|^{-4}X_\flat\wedge{d}X_\flat\otimes\nabla_{*\nabla{}X}X, \\
&\quad{}Tr\left(K\wedge{D}K\right)
\;=\;-\frac{\left<\nabla{X},*\nabla{X}\right>}{|X|^4}X_\flat\wedge{d}X_\flat \\
&\quad{}*Tr(K\wedge*{F})
\;=\;|X|^{-2}X_\flat\wedge\left<*F,\nabla{X}\right> \;=\;|X|^{-2}X_\flat\wedge\left<W(\cdot,\,\cdot),\nabla{X}\right>.
\end{aligned} \label{EqnsComponentsOfTransgression}
\ee
where $X_\flat\wedge\left<W(\cdot,\cdot),\nabla{}X_\flat\right>$ means $X^sg_{si}W_{jkp}{}^qX^p{}_{,q}dx^i\wedge{}dx^j\wedge{}dx^k$.
Letting $\mathcal{P}=\mathcal{P}_\chi$ be the Euler 4-form, the symmetric polynomial is $\mathcal{P}_\chi(h_1,h_2)=(8\pi^2)^{-1}Tr(h_1*h_2)$ for $h_1,h_2\in\mathfrak{so}(4)$.
From (\ref{EqnTPDim4}) and (\ref{EqnsComponentsOfTransgression}) we arrive at
\begin{eqnarray}
	\quad\quad \mathcal{TP}_\chi
	=\frac{1}{8\pi^2}\left(2|X|^{-2}X_\flat\wedge\left<W(\cdot,\cdot),*\nabla{X}\right>
	+|X|^{-4}\left<\nabla{X},*\nabla{X}\right>\,X_\flat\wedge{d}X_\flat\right).
	\label{EqnEulerTransgression}
\end{eqnarray}
Obviously this is singular only on the null-set of $X$.

\subsection{Integrating the transgression when $|\nabla{X}|/|X|$ is large} \label{SubSecFlatCase}

Consider the space of Killing fields $|X|$ in flat $\mathbb{R}^4$, for which $|\nabla{X}|/|X|$ reaches an extremum at the origin.
The tensor fields $X$ and $d{X}_\flat$ can be expressed in a natural basis as
\beg
\begin{aligned}
		X&=\epsilon_1\left(x^1\frac{\partial}{\partial{}x^2}-x^2\frac{\partial}{\partial{}x^1}+x^3\frac{\partial}{\partial{}x^4}-x^4\frac{\partial}{\partial{}x^3}\right)
		+\delta_1\left(x^1\frac{\partial}{\partial{}x^2}-x^2\frac{\partial}{\partial{}x^1}-x^3\frac{\partial}{\partial{}x^4}+x^4\frac{\partial}{\partial{}x^3}\right) \\
		&\;+\epsilon_2\left(x^1\frac{\partial}{\partial{}x^3}-x^3\frac{\partial}{\partial{}x^1}-x^2\frac{\partial}{\partial{}x^4}+x^4\frac{\partial}{\partial{}x^2}\right)
		+\delta_2\left(x^1\frac{\partial}{\partial{}x^3}-x^3\frac{\partial}{\partial{}x^1}+x^2\frac{\partial}{\partial{}x^4}-x^4\frac{\partial}{\partial{}x^2}\right) \\
		&\;+\epsilon_3\left(x^1\frac{\partial}{\partial{}x^4}-x^4\frac{\partial}{\partial{}x^1}+x^2\frac{\partial}{\partial{}x^3}-x^3\frac{\partial}{\partial{}x^2}\right) +\delta_3\left(x^1\frac{\partial}{\partial{}x^4}-x^4\frac{\partial}{\partial{}x^1}-x^2\frac{\partial}{\partial{}x^3}+x^3\frac{\partial}{\partial{}x^2}\right) \\
		&\;+\gamma_1\frac{\partial}{\partial{}x^1}+\gamma_2\frac{\partial}{\partial{}x^2}+\gamma_3\frac{\partial}{\partial{}x^3}+\gamma_4\frac{\partial}{\partial{}x^4} \\
		\frac12dX_\flat
		&=\,\epsilon_1\left(dx^1\wedge{}dx^2+dx^3\wedge{}dx^4\right)	+\delta_1\left(dx^1\wedge{}dx^2-dx^3\wedge{}dx^4\right) \\
		&\;+\epsilon_2\left(dx^1\wedge{}dx^3-dx^2\wedge{}dx^4\right)
		+\delta_2\left(dx^1\wedge{}dx^3+dx^2\wedge{}dx^4\right) \\
		&\;+\epsilon_3\left(dx^1\wedge{}dx^4+dx^2\wedge{}dx^3\right)
		+\delta_3\left(dx^1\wedge{}dx^4-dx^2\wedge{}dx^3\right),
	\end{aligned}
\ees
where the $\epsilon_i$, $\delta_i$ and $\gamma_i$ are constants.
The norms are:
\begin{eqnarray}
	\begin{aligned}
		&|d{X}_\flat|^2=4\left(\boldsymbol\epsilon^2+\boldsymbol\delta^2\right) \\
		&|X_\flat|^2=\left(\boldsymbol\epsilon^2+\boldsymbol\delta^2\right)r^2+\boldsymbol\gamma^2 \\
		&\left<d{X}_\flat,*d{X}_\flat\right>=4\left(\boldsymbol\epsilon^2-\boldsymbol\delta^2\right)
	\end{aligned} \label{EqnsNorms}
\end{eqnarray}
where $\boldsymbol\epsilon=\sqrt{\epsilon_1{}^2+\epsilon_2{}^2+\epsilon_3{}^2}$, $\boldsymbol\delta=\sqrt{\delta_1{}^2+\delta_2{}^2+\delta_3{}^2}$, and $\boldsymbol\gamma=\sqrt{\gamma_1{}^2+\gamma_2{}^2+\gamma_3{}^2+\gamma_4{}^2}$ are constants, and $r=\sqrt{(x^1)^2+(x^2)^2+(x^3)^2+(x^4)^2}$ is the radial function.
A calculation shows that
\begin{eqnarray}
	\begin{aligned}
		&X_\flat\wedge{}dX_\flat=2\left(\boldsymbol\epsilon^2-\boldsymbol\delta^2\right)\left[\left(x^1-G^1\right)dx^2\wedge{}dx^3\wedge{}dx^4 \right.\\ 
		&\quad\quad\quad\quad\quad\quad-\left(x^2-G^2\right)dx^1\wedge{}dx^3\wedge{}dx^4
		+\left(x^3-G^3\right)dx^1\wedge{}dx^2\wedge{}dx^4 \\
		&\quad\quad\quad\quad\quad\quad\left.-\left(x^4-G^4\right)dx^1\wedge{}dx^2\wedge{}dx^3\right]
	\end{aligned}
\end{eqnarray}
where the constants $G^1$, $G^2$, $G^3$, and $G^4$ are
\begin{eqnarray*}
	\begin{aligned}
		&G^1\;=\;\left(-\gamma_1\left(\epsilon_3-\delta_3\right)+\gamma_2\left(\epsilon_2+\delta_2\right)-\gamma_3\left(\epsilon_1+\delta_1\right)\right)/\left(2(\boldsymbol\epsilon^2-\boldsymbol\delta^2)\right), \\
		&G^2\;=\;\left(-\gamma_1\left(\epsilon_2-\delta_2\right)-\gamma_2\left(\epsilon_3+\delta_3\right)+\gamma_4\left(\epsilon_1+\delta_1\right)\right)/\left(2(\boldsymbol\epsilon^2-\boldsymbol\delta^2)\right), \\
		&G^3\;=\;\left(-\gamma_1\left(\epsilon_1-\delta_1\right)+\gamma_3\left(\epsilon_3+\delta_3\right)-\gamma_4\left(\epsilon_2+\delta_2\right)\right)/\left(2(\boldsymbol\epsilon^2-\boldsymbol\delta^2)\right), \\
		&G^4\;=\;\left(\gamma_2\left(\epsilon_1-\delta_1\right)-\gamma_3\left(\epsilon_2-\delta_2\right)+\gamma_4\left(\epsilon_3-\delta_3\right)\right)/\left(2(\boldsymbol\epsilon^2-\boldsymbol\delta^2)\right).
	\end{aligned}
\end{eqnarray*}
Notice
\begin{eqnarray}
	|G^i|\le\frac32\frac{\boldsymbol\gamma\sqrt{\boldsymbol\epsilon^2+\boldsymbol\delta^2}}{|\boldsymbol\epsilon^2-\boldsymbol\delta^2|}. \label{IneqEstG}
\end{eqnarray}
This calculation can be more compactly abbreviated
\begin{eqnarray}
	\begin{aligned}
	&X_\flat\wedge{}dX_\flat\;=\;2\left(\boldsymbol\epsilon^2-\boldsymbol\delta^2\right)\;*\left(\left(\vec{\bf{r}}-\vec{\bf{G}}\right)\cdot{}d\vec{\bf{r}}\right)
	\end{aligned} \label{EqnCompactPartOfTransgression}
\end{eqnarray}
where $\vec{\bf{r}}=\left(x^1,x^2,x^3,x^4\right)$, $\vec{\bf{G}}=\left(G^1,G^2,G^3,G^4\right),$ $d\vec{\bf{r}}=(dx^1,dx^2,dx^3,dx^4)$, and $*$ is the Euclidean hodge star.
Finally, notice that $max|\nabla{X}|/|X|=4(\boldsymbol\epsilon^2+\boldsymbol\delta^2)/\boldsymbol\gamma^2.$

\begin{lemma}[Euclidean Transgression Integration] \label{LemmaFlatTrans}
	Assume $\mathbb{R}^4$ is flat 4-space with a no-where zero Killing field $X$ with $|\nabla{X}|/|X|$ reaching a maximum at the origin.
	Assume $W\in\bigwedge^2\otimes\mathfrak{so}(4)$ is any $\mathfrak{so}(4)$-valued 2-form.
	Consider the 3-form
	\begin{eqnarray}
		\mathcal{TP}
		=2|X|^{-2}X_\flat\wedge\left<W(\cdot,\cdot),*\nabla{X}\right>+|X|^4\left<\nabla{X},*\nabla{X}\right>X_\flat\wedge{}dX_\flat. \label{EqnThreeForm}
	\end{eqnarray}
	Let $\varphi\in{}C^1(\overline{B_r})$ where $B_r$ is the disk of radius $r$ about the origin, let $\|\frac{\partial\varphi}{\partial\vec{x}}\|$ be the manimum of $\left|\frac{\partial\varphi}{\partial{x}^i}\right|$ over the index $i$ and over all points in $B_r$, let $\|W\|$ be the supremum of the norm of $W$ over $B_r$, and denote $\Gamma=max|\nabla{X}|/|X|$.
	Then
	\begin{eqnarray*}
		\left|\int_{B_r}d\varphi\wedge{}\mathcal{TP}\right|\;\le\;C\left\|\frac{\partial\varphi}{\partial\vec{x}}\right\|\left(\|W\|\left(r^2+\Gamma^{-1}\log\Gamma\right)+r+\Gamma^{-1/2}\log(r\Gamma+1)+\Gamma^{-\frac12}\right) \label{IneqEuclTransEst}
	\end{eqnarray*}
	where $C$ is a constant independent of $X$, $|W|$, and $r$.
\end{lemma}
\begin{proof}
We estimate the two terms of $d\varphi\wedge\mathcal{TP}$ in order.
With $|X|^{-2}d\varphi\wedge{}X_\flat\wedge\left<W(\cdot,\cdot),*\nabla{X}\right>$ and using $|\nabla{X}|/|X|=\frac{4(\boldsymbol\epsilon^2+\boldsymbol\delta^2)}{(\boldsymbol\epsilon^2+\boldsymbol\delta^2)r^2+\boldsymbol\gamma^2}$, we estimate
\begin{eqnarray*}
	\begin{aligned}
		&\left|\int_{B_r}|X|^{-2}d\varphi\wedge{}X_\flat\wedge\left<W(\cdot,\cdot),*\nabla{X}\right>\right| \;\le\;\left\|\frac{\partial\varphi}{\partial\vec{x}}\right\|\sup_{B_r}|W|\int_{B_r}\frac{4(\boldsymbol\epsilon^2+\boldsymbol\delta^2)}{(\boldsymbol\epsilon^2+\boldsymbol\delta^2)r^2+\boldsymbol\gamma^2}dVol_{Eucl.} \\
		&\quad\;=\;4|\mathbb{S}^3|\left\|\frac{\partial\varphi}{\partial\vec{x}}\right\|\|W\|\int_0^r\frac{r^3}{r^2+\frac{\boldsymbol\gamma^2}{\boldsymbol\epsilon^2+\boldsymbol\delta^2}}\,dr \\
		&\quad \;=\;2|\mathbb{S}^3|\left\|\frac{\partial\varphi}{\partial\vec{x}}\right\|\|W\|\left(
		r^2+\frac{\boldsymbol\gamma^2}{\boldsymbol\epsilon^2+\boldsymbol\delta^2}\log\left(
		\frac{\frac{\boldsymbol\gamma^2}{\boldsymbol\epsilon^2+\boldsymbol\delta^2}}{r^2+\frac{\boldsymbol\gamma^2}{\boldsymbol\epsilon^2+\boldsymbol\delta^2}}\right)\right).
	\end{aligned}
\end{eqnarray*}
Because $\Gamma=max|\nabla{X}|/|X|=4\frac{\boldsymbol\epsilon^2+\boldsymbol\delta^2}{\boldsymbol\gamma^2}$, we can estimate this by
\begin{eqnarray}
C\left\|\frac{\partial\varphi}{\partial\vec{x}}\right\|\|W\|\left(
r^2+\Gamma^{-1}\log\Gamma\right).
\end{eqnarray}

Now consider the second term, which from (\ref{EqnCompactPartOfTransgression}) is
\begin{eqnarray*}
	\begin{aligned}
		&2\left|\int_{B_r}|X|^{-4}\left<\nabla{X},\,*\nabla{X}\right>\left(\boldsymbol\epsilon^2-\boldsymbol\delta^2\right)\;\frac{\partial\varphi}{\partial{r}}dr\wedge*\left(\left(\vec{\bf{r}}-\vec{\bf{G}}\right)\cdot{}d\vec{\bf{r}}\right)\right|.
	\end{aligned}
\end{eqnarray*}
Using the estimates (\ref{EqnsNorms}), this is bounded by
\begin{eqnarray*}
	\begin{aligned}
		&8\left\|\frac{\partial\varphi}{\partial\vec{x}}\right\|\int_{B_r}
		\frac{(\boldsymbol\epsilon^2-\boldsymbol\delta^2)^2}{\left((\boldsymbol\epsilon^2+\boldsymbol\delta^2)r^2+\boldsymbol\gamma^2\right)^2} \left(r+|G|\right)\,dVol_{Eucl},
	\end{aligned}
\end{eqnarray*}
and we evaluate this integral to get
\begin{eqnarray}
	\begin{aligned}
		&4|\mathbb{S}^3|\left\|\frac{\partial\varphi}{\partial\vec{x}}\right\|
		\left(\frac{\boldsymbol\epsilon^2-\boldsymbol\delta^2}{\boldsymbol\epsilon^2+\boldsymbol\delta^2}\right)^2
		\int_0^r
		\frac{(r+|G|)\,r^3}{\left(r^2+\frac{\boldsymbol\gamma^2}{\boldsymbol\epsilon^2+\boldsymbol\delta^2}\right)^2}\,dr \\
		&\quad\quad\;=\;
		\left[
		\frac{\boldsymbol\gamma^2}{\boldsymbol\epsilon^2+\boldsymbol\delta^2}\frac{r+|G|}{r^2+\frac{\boldsymbol\gamma^2}{\boldsymbol\epsilon^2+\boldsymbol\delta^2}}
		-|G|
		+|G|\log\left(\frac{r^2+\frac{\boldsymbol\gamma^2}{\boldsymbol\epsilon^2+\boldsymbol\delta^2}}{\frac{\boldsymbol\gamma^2}{\boldsymbol\epsilon^2+\boldsymbol\delta^2}}\right)\right. \\
		&\quad\quad\quad\quad\quad\quad\left.-3\left(\frac{\boldsymbol\gamma^2}{\boldsymbol\epsilon^2+\boldsymbol\delta^2}\right)^{\frac12}\tan^{-1}\left(r/\left(\frac{\boldsymbol\gamma^2}{\boldsymbol\epsilon^2+\boldsymbol\delta^2}\right)^{\frac12}\right)
		+2r
		\right]
	\end{aligned}
\end{eqnarray}
Using $\Gamma=4\frac{\boldsymbol\epsilon^2+\boldsymbol\delta^2}{\boldsymbol\gamma^2}$ and $|G|\frac{|\boldsymbol\epsilon^2-\boldsymbol\delta^2|}{\boldsymbol\epsilon^2+\boldsymbol\delta^2}<\frac32\Gamma^{-\frac12}$ from (\ref{IneqEstG}), we bound this above by
\begin{eqnarray}
	\begin{aligned}
		&C\left\|\frac{\partial\varphi}{\partial\vec{x}}\right\|\left(\Gamma^{-\frac12}+\Gamma^{-\frac12}\log\left(\frac14r^2\Gamma+1\right)+2r\right)
	\end{aligned}
\end{eqnarray}
and the estimate (\ref{IneqEuclTransEst}) follows.
\end{proof}

To be usable, we require a version of Lemma \ref{LemmaFlatTrans} for manifolds that are not flat but are very close to being flat.
Further, we want to use this in the case that $X$ is not necessarily known to have the form that it has in the flat case, but where it is known only that $|\nabla{X}|/|X|$ is very large at a point.
Lemma \ref{LemmaStructOfX} addresses the second issue, and Lemma \ref{LemmaBadBallIntegralEstimate} shows that if $|\nabla{X}|/|X|$ at $p$ is large compared to $r_{\mathcal{R}}(p)$, then we can estimate $\int{}d\varphi\wedge\mathcal{TP}_\chi$ over a ball whose radius is small compared $r_{\mathcal{R}}(p)$.

\begin{lemma} \label{LemmaStructOfX}
	Let $X$ be a Killing field on a Riemannian manifold $(N^4,g)$, and assume there is a constant $A>3$ so that $\frac{|\nabla{X}|_p}{|X|_p}>Ar_{\mathcal{R}}(p)^{-1}$ at some point $p$.
	Then $|X|$ has a local minimum within the ball $B(x,6A^{-1}r_{\mathcal{R}}(p))$.
	Further, within the ball $B(x,\frac12r_{\mathcal{R}}(p))$, the minimum of $|X|$ is obtained precisely on a single connected totally geodesic submanifold of dimension 0 or 2.
\end{lemma}
\begin{proof}
	The quantity $r_{\mathcal{R}}(p)\cdot|\nabla{X}|_p/|X|_p$ is scale-invariant under scaling either the metric or $X$, so we may assume $r_{\mathcal{R}}(p)=1$ and $|X|_p=1$.
	Let $\gamma:[-1,1]\rightarrow{}N^4$ be a unit-speed geodesic with $\gamma(0)=0$ (along $\gamma$ we use $X_t$ to denote the vector $X$ at $\gamma(t)$).
	
	Let $Y$ be the Jacobi field along $\gamma$ with $Y_0=X_0$ and $(\nabla_{\dot\gamma}Y)_0=0$.
	Then the fact that $|\Riem|\le1$ forces $|Y_t|\le\cosh(t)$.
	In particular, on $\partial{B}(p,3A^{-1})$ we have $|Y|\le\cosh(3A^{-1})$.
	
	The field $X-Y$ is Jacobi along $\gamma$, and $|\nabla_{\dot\gamma}(X-Y)_0|>A$ and $(X-Y)_0=0$.
	Rauch comparison, along with $|\Riem|\le1$ along $\gamma$, shows that $|(X-Y)_t|>A|\sin(t)|$ for $t\in[-1,1]$.
	Thus on the sphere $\partial{}B(p,kA^{-1})$ we have $|X|+|Y|\;\ge\;|X-Y|\ge{}A\sin(3A^{-1})$.
	Using a few terms of the Taylor series and setting $t=3A^{-1}$, we have
	\beg
	&&|X|
	\;\ge\;
	A\left(3A^{-1}\,+\,\frac13{}27A^{-3}\right)\,-\,1
	\;=\;\left(2\,+\,9A^{-2}\right)
	\;>\;1, \quad \text{when}\;A>3.
	\ees
	Thus on $\partial{}B(p,3A^{-1})$ we have $|X|>1$, and since $|X|=1$ at $p$, necessarily $|X|$ has a minimum in the interior of the disk $B(p,3A^{-1})$.
	This verifies the first claim.
	
	For the second claim, assume $|X|$ obtains a minimum at two points $q_0,q_1\in{}B(p,1/2)$.
	Let $\gamma$ be a minimal geodesic joining $q_0$, $q_1$; because $|\Riem|<1$ in $B(p,1)$, we may assume $\gamma$ remains within $B(p,1)$ (by the Klingenberg Lemma).
	But since $|X|^2$ has a local extreme $q_0$ and $q_1$, ${\dot\gamma}|X|^2=0$ at both points.
	But $X$ is a Jacobi field, so unless $|X|$ is constant along $\gamma$, necessarily the variation form along $\gamma$ has a critical point between $q_0$ and $q_1$.
	But this is impossible by the usual Index Theorem.
	Thus $|X|$ is constant along $\gamma$, which proves that the locus of local minima of $|X|$ is connected and totally geodesic.
	The dimension claim is standard.
\end{proof}

Now we can estimate $\int_{B_i}d\varphi\wedge\mathcal{TP}_\chi$ in regions where $|\nabla{X}|/|X|$ is large, in an arbitrary Riemannian metric.
\begin{lemma}[Almost-Euclidean Transgression Integration] \label{LemmaBadBallIntegralEstimate}
	Let $(N^4,g)$ be a Riemannian manifold, and let $\varphi$ be a differentiable function.
	Assume $A$ is sufficiently large and $|\nabla{X}|_p/|X|_p\ge{}Ar_{\mathcal{R}}(p)^{-1}$ at a point $p$.
	Then
	\beg
	\left|\int_Bd\varphi\wedge\mathcal{TP}_\chi\right|
	\;\le\;
	C\cdot{}A^{3}\cdot\max_B|d\varphi|\cdot\int_B{}\left(r_{\mathcal{R}}\right)^{-3}\,dVol
	\ees
	where $B=B(p,A^{-1}r_{\mathcal{R}}(p))$ and the constant $C$ is dimensional.
\end{lemma}
The crucial point is that $C$ does not depend on $(N^4,g)$.
\begin{proof}
	Let $B=B(p,A^{-1}r_{\mathcal{R}}(p))$.
	The integral $\int_{B}d\varphi\wedge\mathcal{TP}_\chi$ is scale-invariant, so after scaling distances by $Ar_{\mathcal{R}}(p)^{-1}$, we examining the integral over the unit disk, and we have $|\Riem|\le{}A^{-2}$ on $B(p,A)$.
	We are therefore estimating
	\beg
	\begin{aligned}
		\left|\int_{B_1}d\varphi\wedge\mathcal{TP}_\chi\right|
		\;=\;\frac{1}{8\pi^2}\left|\int_{B_1}2|X|^{-2}X_\flat\wedge\left<W(\cdot,\,\cdot),\,*\nabla{X}\right>
		+|X|^{-4}\left<\nabla{X},\,*\nabla{X}\right>\,X_\flat\wedge{d}X_\flat
		\right|
	\end{aligned}
	\ees
	where $B_1=B(p,1)$ is the unit ball.
	We shall use the notation $B_r$ for $B(p,r)\subset{M}$, and $B(r)=B(o,r)\subset{}T_pM$ for balls centered at the origin in $T_pM$.

	It is possible that $B_1$ has non-trivial topology.
	But since $|\Riem|<A^{-2}$ on the large ball $B_A$, the exponential map is a local homeomorphism $B(A)$ to $B_A$, so we can pull the metric back to a set $K\subset{}T_pM$ with the properties
	\begin{itemize}
		\item[-] $\exp_p:K\rightarrow{}B_p(1)$ is a local homeomorphism that is  $k$-to-$1$
		\item[-] The metric balls $B(1),B(2)\subset{}T_pM$ have the property $B(1)\subset{K}\subset{}B(2)$.
	\end{itemize}
	This set $K$ can be constructed as follows: let $B'\subset{B}(p,1)$ be the ball $B(p,1)$ with the cut locus removed; then there is a uniquely defined set $K_0\subset{}T_pM$ that contains the origin $o\in{}T_pM$, is star-shaped with respect to $o$, and $\exp_p:K_0\rightarrow{}B'$ is a homeomorphism.
	Setting $p_0=0$, let $p_0,p_1,\dots,p_N\in{}B(1)$ be the collection of all points so that $p_i\in{}B(1)$ and $p=\exp_p(p_i)$.
	Let $K_1,K_2,\dots,K_N\subset{}T_pM$ be a collection of all sets so $\exp_p:K_i\rightarrow{}B'$ is a bijective isometry, so $K_i\cap{}K_j=\varnothing$ for $i\ne{j}$, and so $p_i\in{}K_i$.
	Let $K=\overline{\bigcup{}K_i}$.
	Clearly $B(1)\subset{}K$.
	To see that $K\subset{}B(2)$, notice that $\dist(p_0,p_i)<1$ and that $K_i$ is within the ball of radius $1$ around $p_i$, so therefore any point of $K_i$ is within a distance of $2$ from $p_0$.
	
	Since $Vol\,B_1=Vol\,K_0$ and $B(1)\subset{}K$, we have the volue ratio estimate
	\beg
		&&VR\,B_1
		\;=\;\frac{Vol\,B_1}{Vol\,B(1)}
		\;\ge\;\frac{Vol\,K_0}{Vol\,K}
	\ees
	and then we compute
	\beg
		\left|\int_{B_1}d\varphi\wedge\mathcal{TP}_\chi\right|
		=\left|\int_{K_0}d\varphi\wedge\mathcal{TP}_\chi\right| =\frac{Vol{}K_0}{Vol{}K}\left|\int_{K}d\varphi\wedge\mathcal{TP}_\chi\right|
		\le{VR\,B_1}\left|\int_{B(2)}d\varphi\wedge\mathcal{TP}_\chi\right|
	\ees
	We can use Lemmas \ref{LemmaFlatTrans} and \ref{LemmaStructOfX} to estimate the rightmost integral, as long as we correct for the small deviation from $B(1)$ being Euclidean.
	Lemma \ref{LemmaFlatTrans} requires $X$ to have two properties that we do not neceesarily have.
	First, $X$ must be Euclidean; but as $A\rightarrow\infty$, we certainly have that $\mathcal{TP}_\chi$ for $|\Riem|<A^{-2}$ converges on the case that $|\Riem|=0$.
	Second, $X$ must have the maximum of $|\nabla{X}|/|X|$ occuring precisely at the origin.
	But Lemma \ref{LemmaStructOfX} says that if $|\nabla{X}|/|X|$ is large at the origin, the maximum occurs very near by (within a distance of $6A^{-1}$).
	
	Therefore after choosing $A$ sufficiently large, we can estimate that $\left|\int_{B(2)}\varphi\wedge\mathcal{TP}_{\chi}\right|$ is no more than, say, twice the Euclidean equivalent.
	Using $r=2$ in Lemma \ref{LemmaBadBallIntegralEstimate},
	\begin{eqnarray*}
		\begin{aligned}
			\left|\int_{B(2)}d\varphi\wedge{}\mathcal{TP}\right| &\;\le\;C\left\|\frac{\partial\varphi}{\partial\vec{x}}\right\|\left(A^{-2}\left(4+A^{-1}\log{}A\right)+4+A^{-1/2}\log(2A+1)+A^{-\frac12}\right) 	\\
			&\;\le\;C\left\|\frac{\partial\varphi}{d\vec{x}}\right\|.
		\end{aligned}
	\end{eqnarray*}
	The final inequality holds since $A^{-1}\log(A)$, $A^{-2}$, etc, are bounded as $A\rightarrow\infty$.
	The coordinate system $(x^1,x^2,x^3,x^4)$ used to evaluate $\|\partial\varphi/\partial\vec{x}\|$ is nearly the Euclidean system on the nearly-euclidean ball $B(2)\subset{}T_pM$.

	Now we un-scale the metric, which is to say we re-multiply distances by $Ar_{\mathcal{R}}(p)^{-1}$.
	The scale-invariant quantity $\|\partial\varphi/\partial\vec{x}\|$ is controlled by the scale-invariant quantity $A^{-1}r_{\mathcal{R}}(p)\|d\varphi\|$ (where $\|d\varphi\|=\sup_{B_{A^{-1}r_{\mathcal{R}}(p)}}|d\varphi|_g$).
	The volume ratio is also scale-invariant, so we obtain
	\begin{eqnarray}
		\begin{aligned}
			&VR\,B(p,A^{-1}r_{\mathcal{R}}(p)){}\int_{B(2A^{-1}r_{\mathcal{R}}(p))}\left|d\varphi\wedge\mathcal{TP}_\chi\right| \\
			&\quad\quad\;\le\;C\cdot\|d\varphi\|\cdot{}VR\,B(p,A^{-1}r_{\mathcal{R}}(p))\cdot{}A^{-1}r_{\mathcal{R}}(p) \\
			&\quad\quad\;=\;C\cdot\|d\varphi\|\cdot{}Vol\,B(p,A^{-1}r_{\mathcal{R}}(p))\cdot{}A^3r_{\mathcal{R}}(p)^{-3}.
		\end{aligned} \label{IneqTransEstCurv}
	\end{eqnarray}
	The final inequality uses that the volume of the ball $B(A^{-1}r_{\mathcal{R}}(p))$ in the tangent space is close to the Euclidean volume, which is $A^{-4}r_{\mathcal{R}}(p)^4$.
	
	Now we compare the right-hand side of (\ref{IneqTransEstCurv}) with $\int_{B}r_{\mathcal{R}}^{-3}dVol$.
	The function $r_{\mathcal{R}}$ varies only slightly on $B_{A^{-1}r_{\mathcal{R}}(p)}$.
	In fact $|r_{\mathcal{R}}^{-1}|\ge{}(1-A^{-1})r_{\mathcal{R}}(p)^{-1}$ on this ball, so
	\beg
	\begin{aligned}
		\int_{B(p,A^{-1}r_{\mathcal{R}}(p))}\,r_{\mathcal{R}}^{-3}dVol
		&\;\ge\;\left(1-A^{-1}\right)r_{\mathcal{R}}(p)^{-3}Vol\,B(p,A^{-1}r_{\mathcal{R}}(p)) \\
		&\;=\;\left(1-A^{-1}\right)A^{-4}r_{\mathcal{R}}(p)\,VR\,B(p,A^{-1}r_{\mathcal{R}}(p)) \\
	\end{aligned}
	\ees
	Therefore
	\beg
	\begin{aligned}
		&\left|\int_{B(p,1)}d\varphi\wedge\mathcal{TP}_\chi\right|
		\;\le\;C\,\frac{A^3}{1-A^{-1}}\|d\varphi\|\int_{B(p,r_{\mathcal{R}}(p))}r_{\mathcal{R}}^{-3}\,dVol.
	\end{aligned}
	\ees
\end{proof}

\section{Proof of Theorems \ref{ThmI} and \ref{ThmIII}}

\subsection{Weak curvature estimate}

From (\ref{EqnEulerTransgression}) $\left|\mathcal{TP}_\chi\right|\le\frac{1}{8\pi^2}\left(2|W|\frac{|\nabla{X}|}{|X|}+\frac12\frac{|\nabla{X}|^3}{|X|^3}\right),$
but we would like $\left|\int{}d\varphi\wedge\mathcal{TP}_\chi\right|\le{}C\sup|d\varphi|\int{}r_{\mathcal{R}}^{-3}dVol$.
We certainly have $|W|<r_{\mathcal{R}}^{-2}$, but there is no {\it pointwise} control over $|\nabla{X}|/|X|$ in terms of $r_{\mathcal{R}}^{-1}$, as can be seen clearly in the flat case.
Rather, we choose some large $A<\infty$ and segregate the points where $|\nabla{X}|/|X|\le{}A{}r_{\mathcal{R}}^{-1}$ from the points where $|\nabla{X}|/|X|>A{}r_{\mathcal{R}}^{-1}$.
The former we call the ``good'' set $\mathcal{G}$ and the latter the ``bad'' set $\mathcal{B}$.
Formally
\beg
\begin{aligned}
&\mathcal{G}\;=\;\left\{p\in{}N^4\,\big|\,|\nabla{X}|_p/|X|_p\le{}Ar_{\mathcal{R}}(p)^{-1}\right\}, \; and \\
&\mathcal{B}\;=\;\left\{p\in{}N^4\,\big|\,|\nabla{X}|_p/|X|_p>Ar_{\mathcal{R}}(p)^{-1}\right\},
\end{aligned}
\ees
and then of course
\begin{eqnarray*}
	\begin{aligned}
		\left|\int{}d\varphi\wedge\mathcal{TP}_\chi\right|
		&\;\le\;A\,\sup|d\varphi|\int_{\mathcal{G}}r_{\mathcal{R}}{}^{-3} \,+\,\left|\int_{\mathcal{B}}d\varphi\wedge\mathcal{TP}_\chi\right|,
	\end{aligned}
\end{eqnarray*}
and we must estimate the second integral.
To do so, we cover $\mathcal{B}$ with balls $\{B_i\}$, where $B_i=B(p_i,A^{-1}r_{\mathcal{R}}(p_i))$ and $p_i\in\mathcal{B}$.
A Gromov-style covering process can ensure that the multiplicity of such a cover is finite---in fact the standard argument goes through with minimal adaptation.
To see this, choose the half-radius balls $B(p_i,\frac12A^{-1}r_{\mathcal{R}}(p))$ maximally with the property that they are non-intersecting.
Then a volume-comparison argument shows that only uniformly finitely may balls $B(p_j,\frac12A^{-1}r_{\mathcal{R}}(p_j))$ can exist within the larger ball $B(p_i,\frac32A^{-1}r_{\mathcal{R}}(p_i))$.
Therefore only finitely many of the balls $B_i$ can possibly intersect at any one point.

Now we estimate
$\left|\int{}d\varphi\wedge\mathcal{TP}_\chi\right|
\;\le\;\left|\int_{\mathcal{G}}d\varphi\wedge\mathcal{TP}_\chi\right|+
\sum_i\left|\int_{B_i}d\varphi\wedge\mathcal{TP}_\chi\right|.
$
By the definition of $\mathcal{B}$ and Lemma \ref{LemmaBadBallIntegralEstimate}, we have
\begin{eqnarray}
	\begin{aligned}
		\left|\int{}d\varphi\wedge\mathcal{TP}_\chi\right|
		&\;\le\;A^{3}\|d\varphi\|\int_{\mathcal{G}}r_{\mathcal{R}}{}^{-3}+
		CA^{3}\|d\varphi\|\sum_i\int_{B_i}r_{\mathcal{R}}{}^{-3}
	\end{aligned}
\end{eqnarray}
and by the finite multiplicity of the cover, we have the required weak estimate
\begin{eqnarray}
	\begin{aligned}
		\int\varphi|W|^2
		\;=\;8\pi^2\int{}d\varphi\wedge\mathcal{TP}_\chi
		\;\le\;C{}A^{3}\|d\varphi\|\int_{supp\,|d\varphi|}r_{\mathcal{R}}{}^{-3}
	\end{aligned} \label{IneqWFirstEst}
\end{eqnarray}
where $\|\varphi\|\triangleq\max_{\supp\,|d\varphi|}|d\varphi|,$ and $C$ is dimensional.

\subsection{The Cheeger-Tian iteration}

In this section we show that $\int_{N^4}|W|^2$ is bounded by a universal constant, irrespective of the Riemannian manifold $(N^4,g)$.

Given a function $f$, its (Hardy-Littlewood) maximal function with cutoff $s$ is
\beg
M_f^s(p)&=&\sup_{0<r<s}\frac{1}{\Vol{B}_p(r)}\int_{B(p,r)}|f|.
\ees
For any subset $\Omega\subset{}N^4$, let its $s$-{\it thickening} be $\Omega^{(s)}\triangleq\{\,p\in{M}\,|\,\dist(p,\Omega)<s\,\}$.
\begin{lemma}\label{LemmaImprovedHolder}
	Assume $\Omega$ is a precompact domain in an $n$-manifold that has $\Ric\ge0$, and assume $0<\alpha<1$.
	There exists a constant $C=C(n,\alpha)$  so that
	\beg
	\left(\frac{1}{|\Omega|}\int_{\Omega}\left(M_g^s\right)^\alpha\;\right)^\frac1\alpha &\le&C\,\frac{1}{|\Omega^{(s)}|}\int_{\Omega^{(s)}}|g|.
	\ees
\end{lemma}
\begin{proof}
	This is the Hardy-Littlewood ``weak (1,1)'' estimate; see \cite{Stein} for example.
\end{proof}
The next step is to establish a link between the curvature scale $r_{\mathcal{R}}$ and the maximal function $M_{|\Riem|^2}$, which is obtained via $\epsilon$-regularity.
\begin{proposition}[Standard $\epsilon$-regularity]\label{PropVREpsReg}
	There exist constants $C=C(n)<\infty$, $\epsilon_0=\epsilon_0(n)>0$ so that if $N^4$ has zero Ricci curvature and
	\beg
	H\;\triangleq\;\frac{1}{r^{-n}\,B(p,r)}\int_{B(p,r)}|\Riem|^{\frac{n}{2}} \;\le\;\epsilon_0,
	\ees
	then $\sup_{B(p,r/2)}|\Riem|\;<\;C\,r^{-2}\,H^\frac2n.$
\end{proposition}
\begin{proof}
	Among the works with this theorem are \cite{Uhl} \cite{Sib} \cite{And1} \cite{TV1} \cite{TV3} \cite{CW}, and \cite{Web1}.
\end{proof}

It is convenient to introduce the {\it s-local energy radius}
\beg
\rho^s(p)&\triangleq&\sup\,\left\{\,0<r<s\;\;\Big|\;\;\frac{1}{r^{-4}\,B(p,r)}\int_{B(p,r)}|\Riem|^2\;\le\;\epsilon_0\;\right\}.
\ees
We use $\rho(p)$ for $\rho^\infty(p)$.
After possibly choosing $\epsilon_0$ to be smaller, Proposition \ref{PropVREpsReg} directly implies $r^s_{\mathcal{R}}(p)\ge\frac12\rho^s(p)$.
A useful way of stating this is the following (see section 4 of \cite{CT}).
\begin{lemma}[Curvature radius weak estimate]\label{LemmaLowerrsRBound}
	Assume $(N^n,g)$ has $\Ric\ge0$.
	Then given any $s>0$, $k\ge0$ we have
	\begin{eqnarray}
	\left(r_{\RRR}(p)\right)^{-k}
	\;\le\;\left(r^s_{\RRR}(p)\right)^{-k}
	\;\le\;\max\left\{2^ks^{-k}\, , \; \left(16\,\epsilon_0^{-1}\,M^s_{|\Riem|^{2}}(p)\right)^{\frac{k}{4}}\right\}. \label{IneqCurvScaleEst}
	\end{eqnarray}
\end{lemma}
\begin{proof}
	Proposition \ref{PropVREpsReg} certainly shows $r^s_{\mathcal{R}}(p)\ge\frac12\rho^s(p)$.
	First if $\rho^s(p)=s$ then $r_{\RRR}^s(p)\ge\frac12s$.
	Then assuming $\rho^s(p)<s$, we have directly from the definitions
	\begin{eqnarray*}
		\epsilon_0\left(\rho^s(p)\right)^{-4} \;=\;\frac{1}{\Vol{B}(p,\rho^s(p))}\int_{B(p,\rho^s(p))}|\Riem|^{2} \;\le\; M^s_{|\Riem|^{2}}(p).\label{IneqBoundingrByM}
	\end{eqnarray*}
	Therefore 
	$\left(r_{\RRR}(p)\right)^{-4}
	\le\left(r^s_{\RRR}(p)\right)^{-4}
	\le16\left(\rho^s(p)\right)^{-4}
	\le16\epsilon_0^{-1}M^s_{|\Riem|^{2}}(p).$
\end{proof}

\begin{lemma} \label{LemmaWeylEst}
	Assume $(N^4,g)$ is a Ricci-flat 4-manifold with a Killing field $X$.
	Let $\Omega$ be a domain with $s$-thickening $\Omega^{(s)}$, and assume $X$ has no zeros in $\Omega^{(s)}$.
	Then
	\begin{eqnarray}
		\int_{\Omega}|W|^2
		&\le&Cs^{-4}|\Omega^{(s)}\setminus\Omega|
		+Cs^{-1}|\Omega^{(s)}\setminus\Omega|^{\frac14}\left(\int_{\Omega^{(s)}\setminus\Omega}|W|^2\right)^{\frac34} \label{IneqEnergyEstFirst} \\
		&\le&Cs^{-4}|\Omega^{(s)}|
		+Cs^{-1}|\Omega^{(s)}|^{\frac14}\left(\int_{\Omega^{(s)}}|W|^2\right)^{\frac34}. \label{IneqEnergyEstSecond}
	\end{eqnarray}
\end{lemma}
\begin{proof}
	Let $\varphi\in{}C^1$ be a cutoff function, with $\varphi\equiv1$ on $\Omega^{(s/4)}$, $\varphi\equiv0$ outside $\Omega^{(3s/4)}$, and $|d\varphi|\le4s^{-1}$.
	Because $\mathcal{TP}_\chi$ is smooth where $|X|\ne0$, it is smooth on $supp\,\varphi$.
	Thus Stokes gives $\int\varphi\mathcal{P}_{\chi}=-\int{}d\varphi\wedge{}\mathcal{TP}_\chi,$ and with $\mathcal{P}_{\chi}=\frac{1}{32\pi^2}|W|^2$, the estimate (\ref{IneqWFirstEst}) gives
	\begin{eqnarray}
		\int_{\Omega}|W|^2
		\;\le\;Cs^{-1}\int_{\Omega^{(3s/4)}\setminus\Omega^{(s/4)}}\left(r_{\mathcal{R}}\right)^{-3}\,dVol
	\end{eqnarray}
	where $C$ is a dimensinal constant.
	Lemma \ref{LemmaLowerrsRBound} gives
	\begin{eqnarray}
		\begin{aligned}
			&s^{-1}\int_{\Omega^{(3s/4)}\setminus\Omega^{(s/4)}}\left(r_{\mathcal{R}}\right)^{-3}\,dVol
			\;\le\;Cs^{-1}\int_{\Omega^{(3s/4)}\setminus\Omega^{(s/4)}}\left(s^{-3}+M^{s/4}_{|W|}{}^{3/4}\right)\,dVol \\
			&\quad\quad\;\le\;Cs^{-4}|\Omega^{(3s/4)}\setminus\Omega^{(s/4)}|
			+Cs^{-1}\frac{|\Omega^{(3s/4)}\setminus\Omega^{(s/4)}|}{|\Omega^{(s)}\setminus\Omega|^{3/4}}\left(\int_{\Omega^{(s)}\setminus\Omega}|W|^2\right)^{\frac34} \\
			&\quad\quad\;\le\;Cs^{-4}|\Omega^{(s)}\setminus\Omega|
			+Cs^{-1}|\Omega^{(s)}\setminus\Omega|^{\frac14}\left(\int_{\Omega^{(s)}\setminus\Omega}|W|^2\right)^{\frac34}
		\end{aligned}
	\end{eqnarray}
\end{proof}

Proceeding with the proof of Theorems \ref{ThmI} and \ref{ThmIII}, choose a radius $R$ so big that the zero locus of $X$ lies entirely within $B_R$.
Choose any $s>0$, and choose $\varphi$ so that $\varphi\equiv1$ on the ball of radius $R+2s$, equals zero outside the ball of radius $2R+2s$, and has $|d\varphi|<2/R$.
The standard Chern-Gauss-Bonnet argument gives $\int\varphi|W|^2
\;=\;\sum_i\chi(N_{X,i})\,-\,\int{}d\varphi\wedge\mathcal{TP}$ and by (\ref{IneqEnergyEstFirst}), we have
\begin{eqnarray}
	\begin{aligned}
		\left|\int{}d\varphi\wedge\mathcal{TP}\right|
		&\le\,\frac{C}{Rs^3}|A_{R+s,2R+3s}| +\frac{C}{R}|A_{R+s,2R+3s}|^{\frac14}\left(\int_{A_{R+s,2R+3s}}|W|^2\right)^{\frac34} \\
		&=\;\,CR^{-4}|A_{2R,5R}|+CR^{-1}|A_{2R,5R}|^{\frac14}\left(\int_{A_{2R,5R}}|W|^2\right)^{\frac34}
	\end{aligned} \label{EqnsOrigW}
\end{eqnarray}
where in the last line we have choosen $s=R$.

Now we set up the iteration process, with the aim to estimate the integral $\int_{A_{2R,5R}}|W|^2$.
Let $\underline{R}_n=2R-R\frac{6}{\pi^2}\sum_{i=1}^n\frac{1}{i^2}$, and $\overline{R}=5R+R\frac{6}{\pi^2}\sum_{i=1}^n\frac{1}{i^2}$; notice that $\underline{R}_n\searrow{}R$ and $\overline{R}_n\nearrow6R$.
First using (\ref{IneqEnergyEstSecond}) and then H\"older's inequality, we have
\begin{eqnarray*}
	\begin{aligned}
		\int_{A_{\underline{R}_{n-1},\overline{R}_{n-1}}}|W|^2
		&\le{}C(R/n^2)^{-4}|A_{\underline{R}_n,\overline{R}_n}|
		+C(R/n^2)^{-1}|A_{\underline{R}_n,\overline{R}_n}|^{\frac14}\left(\int_{A_{\underline{R}_n,\overline{R}_n}}|W|^2\right)^{\frac34} \\
		&\le{}\left(C+\frac14C^4\right)(R/n^2)^{-4}|A_{R,6R}|
		+\frac34\int_{A_{\underline{R}_n,\overline{R}_n}}|W|^2.
	\end{aligned}
\end{eqnarray*}
Iterating $k$ times gives
\begin{eqnarray*}
	\begin{aligned}
		\int_{A_{\underline{R}_{1},\overline{R}_{1}}}|W|^2
		&\le{}
		\left(C'\sum_{n=2}^k\frac{3^n}{4^n}(n+2)^8\right)R^{-4}|A_{R,6R}|
		+\left(\frac34\right)^{k+1}\int_{A_{\underline{R}_{k},\overline{R}_{k}}}|W|^2
	\end{aligned}
\end{eqnarray*}
where $C'=C+\frac{C^4}{4}$ is dimensional.
But $\sum_{n=1}^\infty(3/4)^nn^8<\infty$, so sending $k\rightarrow\infty$,
\begin{eqnarray*}
	\begin{aligned}
		\int_{A_{\underline{R}_{1},\overline{R}_{1}}}|W|^2
		&\le{}C''R^{-4}|A_{R,6R}|.
	\end{aligned}
\end{eqnarray*}
Combining with (\ref{EqnsOrigW}) gives, for a $C''$ that is independent of $R$,
\begin{eqnarray}
	\begin{aligned}
		\int_{B_R}|W|^2
		&\le{}\sum_i\chi(N_X,i)+C''R^{-4}|A_{R,6R}|.
	\end{aligned}
\end{eqnarray}
Because $R^{-4}|A_{R,6R}|<R^{-4}|B_R|$, and $R^{-4}|B_R|<2\pi^2$ by Bishop-Gromov volume comparison, we see in particular $\int_{N^4}|W|^2<\sum_i\chi(N_X,i)+2\pi^2C''$.
Additionally if $\lim_{R\rightarrow\infty}R^{-4}|B_R|=0$ and the Killing field is nowhere-zero, then the manifold is flat.
Now Theorems \ref{ThmI} and \ref{ThmIII} are verified in the case that $\lim_{R\rightarrow\infty}R^{-4}|B_R|=0$.

\subsection{Consideration of the ALE case} \label{SubSecALECase}

We examine when $\lim_{R\rightarrow\infty}R^{-4}|B_R|>0$.
From above, $\int_{N^4}|W|^2<2\pi^2C''$ where $C''$ is a constant independent of $(N^4,g)$, so the Ricci-flat manifilold $(N^4,g)$ is an ALE manifold with one end \cite{BKN}, \cite{And1}, \cite{And4}.
Passing to the universal cover, we have a simply-connected, one-ended ALE manifold with finite energy.

Consider the case that the Killing field $X$ has no zeros.
First, assume $X$ has bounded orbits; then we must have $\chi(N^4)=0$ by the Lefschetz fixed point theorem.
Letting $b_i(N^4)$ be the Betti numbers of $N^4$, we therefore have $1+b_2(N^4)-b_3(N^4)=0$.
Since $N^4$ is ALE, it is homotopy equivalent to a compact manifold-with-boundary $\overline{N}{}^4$ with a single boundary component $\partial\overline{N}{}^4$, which is topologically a quotient of a 3-sphere.
The relevant portion of the relative cohomology sequence is
\begin{eqnarray*}
\longrightarrow{}H^3(\overline{N},\,\partial\overline{N})
\longrightarrow{}H^3(\overline{N})
\longrightarrow{}H^3(\partial\overline{N})
\longrightarrow{}H^4(\overline{N},\partial\overline{N})
\longrightarrow{}0
\end{eqnarray*}
By Poincare duality, $H^3(\overline{N},\partial\overline{N})=H_1(\overline{N})=0$.
Since $\partial\overline{N}$ has one component, $H^3(\partial{}\overline{N})
\longrightarrow{}H^4(\overline{N},\partial\overline{N})$ is an isomorphism.
This forces $H^3(\overline{N})=0$.
Therefore $b_3(\overline{N})=b_3(N)=0$, and so $\chi(N^4)=1+b_2(N)\ge1$.
This means $X$ cannot have bounded orbits.
Since $\chi(N)\ne0$, $X$ must have unbounded orbits, and therefore $\int|W|^2$ is zero or infinite.
But $\int|W|^2$ is uniformly bounded, and so $\int|W|^2=0$.
This completes the proof of Theorem \ref{ThmI}.

In the case that $X$ has zeros, (\ref{IneqMainRiemEst2}) can be found in \cite{EGH} or \cite{GPR}.

\end{document}